\documentclass[11pt]{article}

\usepackage{amsmath,amssymb,amsthm,graphicx}
\usepackage{enumitem}
\usepackage[dvipsnames]{xcolor}
\usepackage[hidelinks]{hyperref}
\usepackage[export]{adjustbox}
\usepackage{wrapfig}

\newtheorem{theorem}{Theorem}

\theoremstyle{definition}

\theoremstyle{remark}

\definecolor{darkgreen}{RGB}{0,128,0}

\newcommand{\Int}{\operatorname{int}}

\usepackage[margin=0.8in]{geometry} 
\usepackage{titling}
\title{A short proof of the universality of the fundamental shadow links}
\author{Yuxuan Yang}

\begin{document}
\maketitle
\vspace{-2em}


\section{Introduction}

Turaev \cite[Chapter~IX]{Turaev} established the shadow theory of 4-manifolds.
Costantino and Thurston \cite{CostantinoThurston2008} developed this theory and proposed the fundamental shadow link, which is an important object in quantum topology, especially in research related to the volume conjecture.

\begin{figure}[htbp]
    \centering
    \begin{minipage}[t]{0.31\textwidth}
        \centering
        \begin{minipage}[c][4.5cm][c]{\linewidth}
            \centering
            \includegraphics[width=\linewidth, height=4.5cm, keepaspectratio]{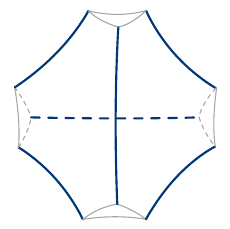}
        \end{minipage}
        \caption{Truncated tetrahedron}
        \label{f1}
    \end{minipage}
    \hfill
    \begin{minipage}[t]{0.31\textwidth}
        \centering
        \begin{minipage}[c][4.5cm][c]{\linewidth}
            \centering
            \includegraphics[width=\linewidth, height=4.5cm, keepaspectratio]{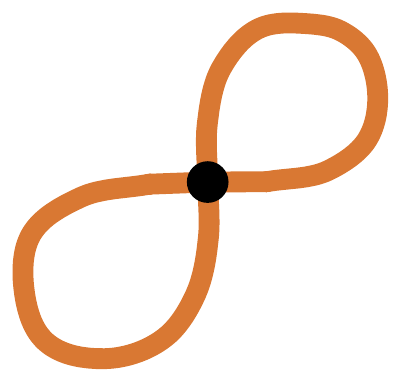}
        \end{minipage}
        \caption{The gluing graph $G$}
        \label{f2}
    \end{minipage}
    \hfill
    \begin{minipage}[t]{0.31\textwidth}
        \centering
        \begin{minipage}[c][4.5cm][c]{\linewidth}
            \centering
            \includegraphics[width=\linewidth, height=4.5cm, keepaspectratio]{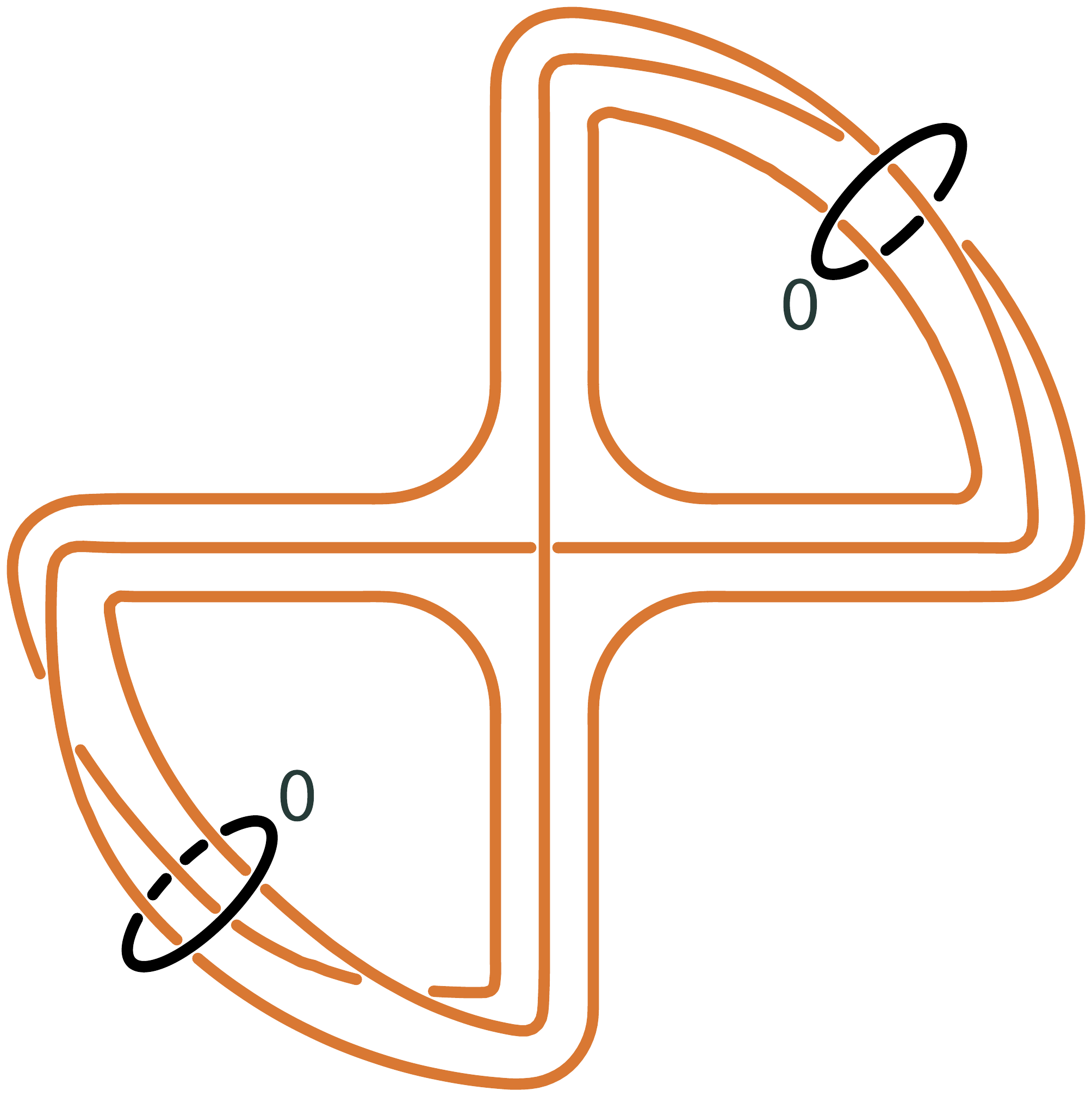}
        \end{minipage}
        \caption{Link diagram of fundamental shadow link}
        \label{f3}
    \end{minipage}
\end{figure}

To define the fundamental shadow link, we start from a truncated tetrahedron $T$ (see Figure \ref{f1}).
Gluing two copies of $T$ along four corresponding hexagonal faces yields a \emph{block}, which is $S^3$ with four open balls removed.
Each boundary sphere of a block is the double of a truncation face.
The six old edges become six properly embedded arcs in the block, with three endpoints on each boundary sphere.

Suppose a connected closed 3-manifold $M_n$ is obtained from $n\geq1$ copies of this block by pairing all boundary spheres using orientation-reversing maps that match the three marked endpoints.
The arcs join together to form the link $L$.
Then we say $(M_n,L)$ is a \emph{fundamental shadow link pair}.

To understand the topology of the ambient manifold $M_n$, let $G$ (see Figure \ref{f2}) be a connected 4-regular graph such that each vertex corresponds to a block in $M_n$ and each edge corresponds to the pairing between boundary spheres of blocks, which we call the gluing graph. $G$ has $n$ vertices and $2n$ edges.
First gluing the blocks along a maximal tree of $G$ produces $S^3$ with $2n+2$ open balls removed.
Each of the remaining $n+1$ sphere pairings adds one $S^2\times S^1$ summand.
Hence $M_n=\#^{n+1}(S^2\times S^1)$.

The following standard steps yield the link diagram of $(M_n,L)$ in $S^3$ (see Figure \ref{f3}).
We can replace each vertex of $G$ with the thick lines in Figure \ref{f1}; then the edges of $G$ should be replaced by pairings of three pairs of points. Recall that the $n+1$ edges not included in the maximal tree of $G$ correspond to the $S^2\times S^1$ summands; we place a circle carrying 0-surgery around each of these edges.
Then these $n+1$ circles form the link diagram of $M_n$, and $G$ has been transformed into a diagram of $L$ (see the orange link in Figure \ref{f3}).

A fundamental result of Costantino and Thurston \cite{CostantinoThurston2008} is the following universality statement for the fundamental shadow links:

\begin{theorem}\label{t1}
Let $M$ be a compact connected oriented 3-manifold whose boundary is empty or a union of tori.
There exist a fundamental shadow link pair $(M_n,L)$ and a partition $L=L_0\sqcup L_1$ such that integral surgery on $L_0$ produces a closed manifold $M'$ satisfying
\[
M\cong M'\setminus\Int\nu L_1,
\]
where $\nu L_1$ is a closed tubular neighbourhood.
\end{theorem}

When $M$ is closed, it is stated in \cite[Proposition~3.36]{CostantinoThurston2008}. Their proof used the shadow theory of 4-manifolds.

We will give a constructive three-dimensional proof, which only requires the properties of Dehn surgery and the Lickorish--Wallace theorem \cite{Lickorish1962,Wallace1960}.

\section{Proof of Theorem \ref{t1}}

By the Lickorish--Wallace theorem, there is a link $L'=\textcolor{red}{L_0'}\sqcup \textcolor{blue}{L_1}$ in $S^3$ such that performing integral Dehn surgery on $\textcolor{red}{L_0'}$ in the complement $S^3\setminus\Int\nu\textcolor{blue}{L_1}$ yields $M$. We can assume $\textcolor{red}{L_0'}\ne\emptyset$, otherwise we can replace $\textcolor{red}{L_0'}=\emptyset$ by an unknot with $1$-surgery. Let $D$ be a link diagram of $L'$.  We may assume $D$ has at least one crossing, otherwise we can perform a Reidemeister~I move. Ignoring the over/under information at each crossing, we obtain a planar 4-regular graph $G$. We can assume $G$ is connected, otherwise we can perform a Reidemeister~II move on $D$. 

\begin{wrapfigure}{r}{0.54\textwidth}
    \centering
    \includegraphics[width=\linewidth]{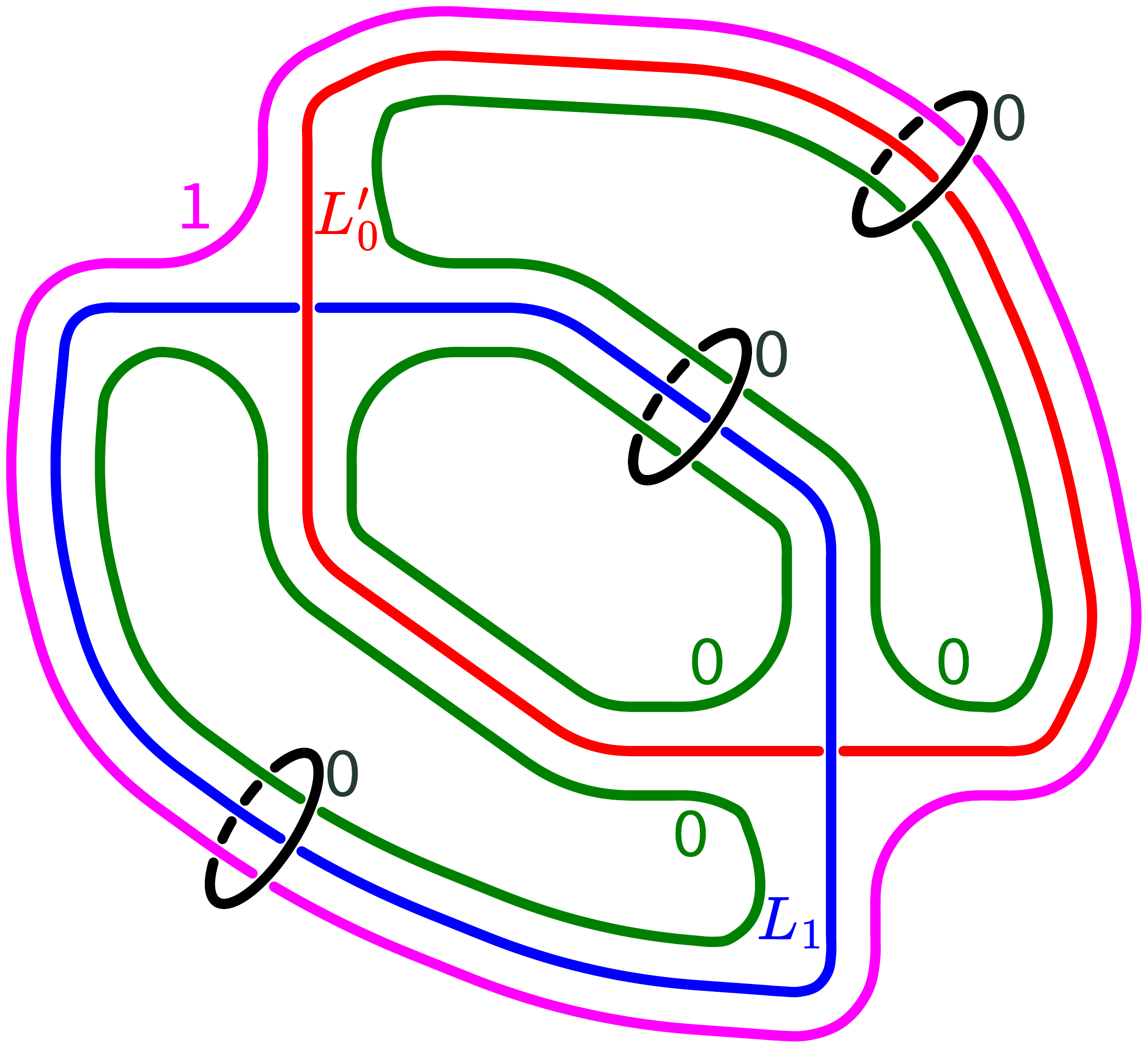}
    \caption{Construct a fundamental shadow link from $D$}
    \label{f5}
\end{wrapfigure}

Suppose $D$ has $n$ crossings. We will construct a fundamental shadow link pair $(M_n,L)$ from $D$.
Let $N$ be a two-dimensional small closed regular neighbourhood of $G$, and consider the link $\partial N$, which consists of one outermost boundary $\textcolor{magenta}{K_{\text{pink}}}$ (see the pink knot in Figure \ref{f5}) and $\textcolor{darkgreen}{L_{\text{green}}}$ with \textcolor{darkgreen}{$n+1$} components corresponding to the boundaries of the interior regions (see the green link in Figure \ref{f5}). 
We equip $\textcolor{magenta}{K_{\text{pink}}}$ with \textcolor{magenta}{$1$-surgery} and each component of $\textcolor{darkgreen}{L_{\text{green}}}$ with \textcolor{darkgreen}{$0$-surgery}. 
Let $L_0=\textcolor{red}{L_0'}\sqcup \textcolor{magenta}{K_{\text{pink}}} \sqcup \textcolor{darkgreen}{L_{\text{green}}}$.
Now we put $n+1$ black circles around the edges which are not included in the maximal tree of $G$ (see the black circles in Figure \ref{f5}) and equip them with 0-surgery. 
Then these $n+1$ circles form the link diagram of $M_n$. 
Let $L=L_0\sqcup \textcolor{blue}{L_1}$. From the standard link diagram (Figure \ref{f5}), we observe that $(M_n,L)$ is a fundamental shadow link pair.

Consider adding $n+1$ edges one by one to the maximal tree of $G$; we find that at each step the edge added (surrounded by a black circle in Figure \ref{f5}) corresponds to a new interior region (bounded by a component of $\textcolor{darkgreen}{L_{\text{green}}}$). 
Thus, as we follow the steps, each black circle corresponds to a component of $\textcolor{darkgreen}{L_{\text{green}}}$. 
Starting from $M_n$ and following the steps of adding edges, we perform \textcolor{darkgreen}{$0$-surgeries} on each component of $\textcolor{darkgreen}{L_{\text{green}}}$, \textcolor{darkgreen}{$n+1$} times in total. 
At each step, a black circle marked with 0 and one component of $\textcolor{darkgreen}{L_{\text{green}}}$ cancel each other. 
So after these \textcolor{darkgreen}{$n+1$} surgeries for $\textcolor{darkgreen}{L_{\text{green}}}$, we are left with $\textcolor{red}{L_0'}\sqcup \textcolor{blue}{L_1}\sqcup \textcolor{magenta}{K_{\text{pink}}}$, which lies in the ambient manifold $S^3$. 
Note that $\textcolor{magenta}{K_{\text{pink}}}$ is unlinked from $\textcolor{red}{L_0'}\sqcup\textcolor{blue}{L_1}$, and that \textcolor{magenta}{$1$-surgery} does not change $S^3$. It follows that $M\cong M'\setminus\Int\nu\textcolor{blue}{L_1}$, where $M'$ is obtained by performing integral surgery on $L_0=\textcolor{darkgreen}{L_{\text{green}}}\sqcup\textcolor{magenta}{K_{\text{pink}}}\sqcup\textcolor{red}{L_0'}$. 
This completes the proof.


\vspace{1.5em}\noindent\textbf{Declaration of AI Use.} The figures in this paper were drawn with the assistance of AI.

\bibliographystyle{unsrt}
\bibliography{references}
\end{document}